\documentclass{amsart}
\usepackage{amssymb}

\newtheorem{theorem}{Theorem}[section]
\newtheorem{proposition}[theorem]{Proposition}
\newtheorem{lemma}[theorem]{Lemma}

\theoremstyle{definition}

\newtheorem{conjecture}[theorem]{Conjecture}

\numberwithin{equation}{section}

\begin{document}
\title[Supercongruences for Truncated Hypergeometric Series]
{Supercongruences for Truncated Hypergeometric Series and $p$-adic Gamma Function}


\author{Rupam Barman}
\address{Department of Mathematics, Indian Institute of Technology Guwahati, North Guwahati, Guwahati-781039, Assam, INDIA}
\curraddr{}
\email{rupam@iitg.ac.in}
\author{Neelam Saikia}
\address{Department of Mathematics, Indian Institute of Technology Guwahati, North Guwahati, Guwahati-781039, Assam, INDIA}
\curraddr{}
\email{neelam16@iitg.ac.in}
\thanks{}


\subjclass[2010]{Primary: 33E50, 33C20; Secondary: 33C99, 11S80.}
\date{To appear in Math. Proc. Cambridge Phil. Soc.}
\keywords{Hypergeometric series, truncated hypergeometric series, $p$-adic Gamma function, supercongruences.}
\thanks{We thank the referee for his/her valuable comments and pointing out mathematical errors in the earlier draft. We thank Ken Ono and Wadim Zudilin for 
commenting on an initial draft of the manuscript. Sincere thanks to Ling Long for many helpful discussions
during the preparation of the article. The first author is supported by a research grant under the MATRICS scheme of SERB, 
Department of Science and Technology, Government of India. The second author acknowledges the financial support of Department of Science and Technology, 
Government of India for supporting a part of this work under INSPIRE Faculty Fellowship.}
\begin{abstract} We prove three more general supercongruences between truncated hypergeometric series and $p$-adic Gamma function
from which some known supercongruences follow. A supercongruence conjectured by Rodriguez-Villegas and
proved by E. Mortenson using the theory of finite field hypergeometric series follows from one of our more general supercongruences.
We also prove a supercongruence for ${_7}F_6$ truncated hypergeometric series
which is similar to a supercongruence proved by L. Long and R. Ramakrishna.
\end{abstract}
\maketitle
\section{Introduction and statement of results}
For a complex number $a$, the rising factorial or the Pochhammer symbol is defined as $(a)_0=1$ and $(a)_k=a(a+1)\cdots (a+k-1), ~k\geq 1$.
If $\Gamma(x)$ denotes the Gamma function, then we have
$(a)_k=\dfrac{\Gamma(a+k)}{\Gamma(a)}$. For a non-negative integer $r$, and $a_i, b_i\in\mathbb{C}$ with $b_i\notin\{\ldots, -3,-2,-1, 0\}$,
the (generalized) hypergeometric series ${_{r+1}}F_{r}$ is defined by
\begin{align}\label{hyper}
{_{r+1}}F_{r}\left[\begin{array}{cccc}
                   a_1, & a_2, & \ldots, & a_{r+1} \\
                    & b_1, & \ldots, & b_r
                 \end{array}; \lambda
\right]:=\sum_{k=0}^{\infty}\frac{(a_1)_k\cdots (a_{r+1})_k}{(b_1)_k\cdots(b_r)_k}\cdot\frac{\lambda^k}{k!},
\end{align}
which converges for $|\lambda|<1$.
\par When we truncate the above sum at $k=n$, it is known as a truncated hypergeometric series. We use subscript notation
to denote the truncated hypergeometric series
\begin{align}
{_{r+1}}F_{r}\left[\begin{array}{cccc}
                   a_1, & a_2, & \ldots, & a_{r+1} \\
                    & b_1, & \ldots, & b_r
                 \end{array}; \lambda
\right]_n:=\sum_{k=0}^{n}\frac{(a_1)_k\cdots (a_{r+1})_k}{(b_1)_k\cdots(b_r)_k}\cdot \frac{\lambda^k}{k!}.\notag
\end{align}
If one of the $a_i$'s is a negative integer, then the hypergeometric series \eqref{hyper} terminates.
More details on hypergeometric series can be found in the books
by W. Bailey \cite{bailey}, L. Slater \cite{slater} and G. Andrews, R. Askey and R. Roy \cite{andrews}.
\par
Hypergeometric series are related to many mathematical objects. For example, they are related to algebraic varieties, differential equations and modular
forms. There are beautiful results expressing periods of abelian varieties in terms of hypergeometric series.
Consider the Legendre family of elliptic curves $E_{\lambda}: y^2=x(x-1)(x-\lambda)$ parameterized by
$\lambda$. A period of $E_{\lambda}$ can be expressed in terms of ${_{2}}F_{1}\left[\begin{array}{cc}
                   \frac{1}{2}, & \frac{1}{2} \\
                    & 1
                 \end{array}; \lambda
\right]$. In \cite{BK}, the first author with G. Kalita defined a period analogue for the algebraic curve
$y^{\ell}=x(x-1)(x-\lambda), \ell \geq 2$;
and described it by a $_{2}F_{1}$-hypergeometric series. In \cite{mccarthy1}, D. McCarthy expressed the real period of the elliptic curve
$y^2=(x-1)(x^2+\lambda)$ by a $_{3}F_{2}$-hypergeometric series. In general, periods are complicated transcendental numbers. However, in
case of CM elliptic curves, the Selberg-Chowla formula predicts that
any period is an algebraic multiple of a quotient of Gamma values \cite{selberg}. In recent years, there has been
a focus on a $p$-adic analog of these complex periods computed from the hypergeometric series and Gamma function. This often leads to
congruences between truncated hypergeometric series and $p$-adic Gamma function modulo $p$. In the case of CM elliptic curves,
stronger congruences have been observed. These congruences are often called as \emph{supercongruences}.
These congruences are stronger than those predicted by commutative formal group theory.
\par In literature, there are many supercongruences conjectured by F. Beukers \cite{beuker2}, van Hamme \cite{vanHamme},
Rodriguez-Villegas \cite{RV}, W. Zudilin \cite{zudilin}, H. Chan et al. \cite{chan}, and Z.-W. Sun \cite{sun1, sun2}.
Interestingly, different methods are needed to prove some of these conjectures. For example, in \cite{ahlgren, ahlgren-ono, mccarthy-fuselier,
kilbourn, mccarthy-osburn, mortenson2, mortenson}, Gaussian hypergeometric series are used to prove some of the supercongruences. The Wilf-Zeilberger method is
used in \cite{zudilin}  and $p$-adic analysis is used in \cite{beuker1, osburn-sahu, swisher}. In \cite{long-pacific}, L. Long proves several supercongruences
related to special valuations of truncated hypergeometric series including a conjecture of van Hamme  using hypergeometric evaluation
identities and combinatorial techniques.
\par
In \cite{RV}, Rodriguez-Villegas conjectured that for any odd prime $p$
\begin{align}\label{cor-MT-5}
{_2}F_{1}\left[\begin{array}{cc}
                          \frac{1}{2}, & \frac{1}{2} \\
                          ~ & 1
                        \end{array};1
\right]_{p-1}:=\sum_{k=0}^{p-1}{2k\choose k}^2\cdot 16^{-k}
\equiv-\Gamma_p\left(\frac{1}{2}\right)^{2}\equiv \left(\frac{-1}{p}\right)\pmod{p^2},
\end{align}
where $\Gamma_p(\cdot)$ denotes the $p$-adic Gamma function recalled in Section 2.
The above congruence was proved by E. Mortenson in \cite{mortenson2}
using the theory of finite field hypergeometric series \cite{greene}. Here, we prove the following
result which generalizes \eqref{cor-MT-5}. We note that our method does not use finite field hypergeometric series.
\begin{theorem}\label{MT-5}
Let $p$ be an odd prime such that $p\equiv1\pmod {n(n-1)}$. Then modulo $p^2$ we have
\begin{align}
{_n}F_{n-1}\left[\begin{array}{cccc}
                          \frac{1}{n}, & \frac{1}{n}, & \ldots, & \frac{1}{n} \vspace{.08cm}\\
                          ~ & \frac{1}{n-1}, & \ldots, & \frac{1}{n-1}
                        \end{array};1
\right]_{\frac{p-1}{n-1}}&\equiv {_n}F_{n-1}\left[\begin{array}{cccc}
                          \frac{1}{n}, & \frac{1}{n}, & \ldots, & \frac{1}{n} \vspace{.08cm}\\
                          ~ & \frac{1}{n-1}, & \ldots, & \frac{1}{n-1}
                        \end{array};1
\right]_{\frac{p-1}{n}}\notag\\
&\equiv(-1)^n\Gamma_p\left(\frac{1}{n-1}\right)^{n-1}\Gamma_p\left(\frac{n-1}{n}\right)^{n}.\notag
\end{align}
\end{theorem}
If we put $n=2$ in Theorem \ref{MT-5}, then \eqref{cor-MT-5} follows.
\par The following supercongruence is obtained by L. Long \cite[Theorem 1]{long-pacific} for any prime $p>3$.
\begin{align}\label{cor-MT-4}
&{_{5}}F_{4}\left[\begin{array}{ccccc}
                 \frac{1}{2}, & \frac{1}{2}, & \frac{1}{2}, & \frac{1}{2}, & \frac{5}{4} \vspace{.08cm}\\
                 ~ & \frac{1}{4}, & 1, & 1, & 1
               \end{array};1
\right]_{p-1}:=\sum_{k=0}^{p-1}(4k+1){2k\choose k}^4\cdot 256^{-k}\equiv p\pmod {p^4}.
\end{align}
We prove the followng supercongruence which generalizes \eqref{cor-MT-4}.
\begin{theorem}\label{MT-4}
Let $p\geq 5$ be a prime such that $p\equiv1\pmod n$. Then modulo $p^4$ we have
\begin{align}
{_5}F_4\left[\begin{array}{ccccc}
               \frac{1}{n}, & \frac{1}{n}, & \frac{1}{n}, & \frac{1}{n}, & 1+\frac{1}{2n} \vspace{.08cm}\\
               ~ & \frac{1}{2n}, & 1, & 1, & 1
             \end{array};1
\right]_{p-1}\equiv
(-1)^{1+\frac{p-1}{n}} p\Gamma_p\left(\frac{1}{n}\right)^2\Gamma_p\left(\frac{n-2}{n}\right).\notag
\end{align}
\end{theorem}
If we put $n=2$ in Theorem \ref{MT-4}, then \eqref{cor-MT-4} follows.
\par We also prove the following supercongruence.
\begin{theorem}\label{MT-6}
Let $p \geq 7$ be a prime. Let $n\in\mathbb{Z}^{+}$ be even such that $n\neq 2, 4$, $p\equiv1\pmod {n}$ and $p\not\equiv1\pmod{2n}$. Then we have
\begin{align}
&_{5}F_{4}\left[\begin{array}{ccccc}
                 \frac{1}{n}, & \frac{1}{n}, & \frac{1}{n}, & \frac{1}{n}, & \frac{2n-3}{2n} \vspace{.08cm}\\
                 ~ & \frac{5}{2n}, & 1, & 1, & 1
               \end{array};1
\right]_{p-1}\notag\\
&\equiv (-1)^{1+\frac{p-1}{n}}\frac{p}{2n} \frac{\Gamma_p\left(\frac{1}{n}\right)^2\Gamma_p\left(\frac{1}{2n}\right)^4
\Gamma_p\left(\frac{2n-3}{2n}\right)^3
\Gamma_p\left(\frac{5}{2n}\right)}{\Gamma_p\left(\frac{2}{n}\right)}\pmod {p^4}.\notag
\end{align}
\end{theorem}
\par Recently, L. Long and R. Ramakrishna \cite{long} prove several supercongruences using a
technique which relies on the relations between the classical and $p$-adic Gamma functions. They also prove a conjecture of J. Kibelbek \cite{kibelbek}
and a strengthened version of a conjecture of van Hamme. For instance, they prove the following supercongruence modulo $p^6$
\cite[Theorem 2]{long} which is stronger than a prediction of van Hamme in \cite{vanHamme}.
 \begin{align}\label{long-ravi}
 \displaystyle {_{7}}F_{6}\left[\begin{array}{ccccccc}
                   \frac{7}{6}, & \frac{1}{3}, &\frac{1}{3}, &\frac{1}{3}, &\frac{1}{3}, &\frac{1}{3}, &\frac{1}{3} \vspace{.08cm}\\
                   &\frac{1}{6}, & 1, &1,&1,&1,&1
                 \end{array}; 1
\right]_{p-1}\equiv \left\{ \begin{array}{ll}
-p\Gamma_p\left(\dfrac{1}{3}\right)^9 \hspace{.6cm}\hbox{if~ $p \equiv 1\hspace{-.25cm}\pmod{6}$;} \\
                              -\dfrac{10}{27}p^4\Gamma_p\left(\dfrac{1}{3}\right)^9    \hbox{if ~$p\equiv 5\hspace{-.25cm}\pmod{6}$.}
                               \end{array}
                             \right.
\end{align}
\par
Similar to the supercongruence \eqref{long-ravi}, we prove the following supercongruence for
${_{7}}F_{6}$ truncated hypergeometric series modulo $p^6$.
\begin{theorem}\label{MT-2}
Let $p\equiv1\pmod 8$ be a prime. Then we have
\begin{align}
{_7}F_6\left[\begin{array}{ccccccc}
               \frac{1}{8}, & \frac{17}{16}, & \frac{1}{4}, & \frac{1}{4},
                & \frac{1}{4}, & \frac{1}{4}, & \frac{1}{4} \vspace{.08cm}\\
               ~ & \frac{1}{16}, & \frac{7}{8}, & \frac{7}{8},
               & \frac{7}{8}, & \frac{7}{8}, & \frac{7}{8}
             \end{array};1
\right]_{\frac{7(p-1)}{8}}\equiv-p\Gamma_p\left(\frac{7}{8}\right)^6\Gamma_p\left(\frac{3}{8}\right)^{10}\pmod {p^6}.\notag
\end{align}
\end{theorem}
\par One of the main results of A. Deines et al. is the following supercongruence \cite[Theorem 7]{DFLST}
which is true for $p\equiv 1 \pmod{4}$.
 \begin{align}\label{fuselier-thm-7}
&_{4}F_{3}\left[\begin{array}{ccccc}
                 \frac{1}{4}, & \frac{1}{4}, & \frac{1}{4}, & \frac{1}{4} \\
                 ~ & 1, & 1, & 1
               \end{array};1
\right]_{p-1}\equiv (-1)^{\frac{p-1}{4}}\Gamma_p\left(\frac{1}{2}\right)\Gamma_p\left(\frac{1}{4}\right)^6\pmod {p^4}.
\end{align}
We now conjecture the following supercongruence which generalizes the supercongruence \eqref{fuselier-thm-7} based on numerical evidence.
\begin{conjecture}\label{conj-1}
 Let $p\geq 7$ be a prime. If $n\in\mathbb{Z}^{+}$ satisfies $p\equiv1\pmod {2n}$ then we have
\begin{align}
&_{5}F_{4}\left[\begin{array}{ccccc}
                 \frac{1}{n}, & \frac{1}{n}, & \frac{1}{n}, & \frac{1}{n}, & \frac{2n-3}{2n} \vspace{.08cm}\\
                 ~ & \frac{5}{2n}, & 1, & 1, & 1
               \end{array};1
\right]_{p-1}\notag\\
&\equiv -\frac{\Gamma_p\left(\frac{1}{n}\right)^2\Gamma_p\left(\frac{1}{2n}\right)^4\Gamma_p\left(\frac{2n-3}{2n}\right)^3
\Gamma_p\left(\frac{5}{2n}\right)}{\Gamma_p\left(\frac{2}{n}\right)}\pmod {p^4}.\notag
\end{align}
\end{conjecture}
Now, if we put $n=4$ in Conjecture \ref{conj-1}, one can deduce \eqref{fuselier-thm-7} for $p\equiv 1 \pmod{8}$ using the product formula
\eqref{prod-formula} for $m=2$ and $x=\frac{1}{4}$. If $n$ is odd, the statement of Conjecture \ref{conj-1} is in fact true for $p\equiv 1 \pmod{n}$.
If $n$ is even, then Theorem \ref{MT-6} includes certain primes which are missed in Conjecture \ref{conj-1}. We remark that
due to the presence of the factor $p$ on the right hand side of Theorem \ref{MT-6} the same technique couldn't be used to prove Conjecture \ref{conj-1}.
\section{Preliminaries}
We first recall the definition of the $p$-adic Gamma function and list some of its main properties.
For further details, see \cite{kob}. Let $\mathbb{Z}_p$ denote the ring of $p$-adic integers and $\mathbb{Q}_p$ denote
the field of $p$-adic numbers. Let $v_p$ and $|\cdot|_p$ denote the $p$-adic valuation and absolute value on $\mathbb{Q}_p$,
respectively.
The $p$-adic gamma function $\Gamma_p$ is defined by setting $\Gamma_p(0)=1$, and for $n \in\mathbb{Z}^+$ by
\begin{align}
\Gamma_p(n):=(-1)^n\prod_{\substack{0<j<n\\p\nmid j}}j.\notag
\end{align}
The function has a unique extension to a continuous function
$\Gamma_p: \mathbb{Z}_p \rightarrow \mathbb{Z}_p^{\times}$. If $x\in \mathbb{Z}_p$ and $x\neq 0$, then
$\Gamma_p(x)$ is defined as
\begin{align}
\Gamma_p(x):=\lim_{x_n\rightarrow x}\Gamma_p(x_n),\notag
\end{align}
where $x_n$ runs through any sequence of positive integers $p$-adically approaching $x$.
This function is locally analytic and has
a Taylor series expansion
$$\Gamma_p(x+z)=\sum_{n=0}^{\infty}a_nz^n, ~~a_n\in \mathbb{Q}_p,$$
with radius of convergence $\varrho=p^{-\frac{1}{p}-\frac{1}{p-1}}$.
\par We now list some of the main properties of $\Gamma_p$ in the following proposition.
\begin{proposition}\label{prop-1}
Let $x\in\mathbb{Z}_p$. We have
\begin{enumerate}
\item $\Gamma_p(0)=1$ and $\Gamma_p(1)=-1$.
\item $ \dfrac{\Gamma_p(1+x)}{\Gamma_p(x)}= \left\{
        \begin{array}{ll}
          -x, & \hbox{if $|x|_p=1$;} \\
          -1, & \hbox{if $|x|_p<1$.}
        \end{array}
      \right. $
\item $\Gamma_p(1-x)\Gamma_p(x)=(-1)^{a_0(x)}$, where $a_0(x)\in\{1,2,\ldots,p\}$ such that $x\equiv a_0(x)\pmod p$.
\item $\Gamma_p\left(\dfrac{1}{2}\right)^2=(-1)^{\frac{p+1}{2}}$.
\end{enumerate}
\end{proposition}
For $x\in \mathbb{Z}_p$ and $n\in \mathbb{N}$, we define
$$G_n(x):=\dfrac{\Gamma_p^{(n)}(x)}{\Gamma_p(x)}.$$
In particular, $G_0(x)=1$. We now state two results about Taylor series of $\Gamma_p$ from \cite{long}. Let $\mathbb{C}_p$ denote the
completion of the algebraic closure of $\mathbb{Q}_p$.
\begin{proposition}\cite[Proposition 13]{long}\label{prop-2}
Let $p\geq 5$ be a prime and $a\in \mathbb{Q}$ with $v_p(a)\geq 0$. Then
$v_p\left(\frac{G_i(a)}{i!}\right)\geq -i\left(\frac{1}{p}+\frac{1}{p-1}\right)$. For $i<p$, $v_p\left(\frac{G_i(a)}{i!}\right)=0$.
We may extend the domain of $\Gamma_p(a+x)$ by setting $\Gamma_p(a+x)=\Gamma_p(a)\cdot \displaystyle \sum_{k=0}^{\infty}\frac{G_k(a)}{k!}x^k$
for $x\in \mathbb{C}_p$ with $v_p(x)\geq \left(\frac{1}{p}+\frac{1}{p-1}\right)$. In particular,
\begin{align}\label{new-eq-1}
\frac{\Gamma_p(a+1+x)}{\Gamma_p(a+x)}=\left\{
         \begin{array}{ll}
           -(a+x), & \hbox{if $|a+x|_p=1$;} \\
           -1, & \hbox{if $|a+x|_p<1$,}
         \end{array}
       \right.
       \end{align}
       and $\Gamma_p(a+x)\Gamma_p(1-a-x)=(-1)^{a_0(a)}$.
\end{proposition}
\begin{theorem}\cite[Theorem 14]{long}\label{thm-long-14}
For $p\geq 5, r\in \mathbb{N}, a\in \mathbb{Z}_p, m\in \mathbb{C}_p$ satisfying $v_p(m)\geq 0$ and $t\in \{0, 1, 2\}$ we have
\begin{align}
 \frac{\Gamma_p(a+mp^r)}{\Gamma_p(a)}\equiv \sum_{k=0}^t\frac{G_k(a)}{k!}(mp^r)^k \pmod{p^{(1+t)r}}.
\end{align}
The above result also holds for $t= 4$ if $p\geq 11$.
\end{theorem}
\begin{lemma}\cite[Lemma 2.1]{swisher}\label{lemma-swisher}
Let $p$ be prime and $\zeta$ a primitive $n$-th root of unity for some positive integer $n$. If $a, b\in \mathbb{Q}\cap \mathbb{Z}_p^{\times}$
and $k$ is a positive integer such that $(a+j)\in \mathbb{Z}_p^{\times}$ for each $0\leq j \leq k-1$, then
$$(a-bp)_k(a-b\zeta p)_k \cdots (a-b\zeta^{n-1}p)_k\equiv (a)_k^n \pmod{p^n}.$$
Moreover for an indeterminate $x$,
$$(a-bx)_k(a-b\zeta x)_k \cdots (a-b\zeta^{n-1}x)_k \in \mathbb{Z}_p[[x^n]].$$
\end{lemma}
Let $x\in \mathbb{C}_p$ be such that $v_p(x)\geq \left(\frac{1}{p}+\frac{1}{p-1}\right)$ and $|x+j|_p=1$
for each $0\leq j\leq k-1$. Then using \eqref{new-eq-1} repeatedly, we deduce that
\begin{align}\label{new-eq-2}
 (x)_k=(-1)^k\frac{\Gamma_p(x+k)}{\Gamma_p(x)}.
\end{align}
If for $j=j_1, j_2, \ldots, j_i$, $|x+j|_p<1$ and $|x+j|_p=1$ for all other values of $j$, then  using \eqref{new-eq-1} repeatedly,
we deduce that
\begin{align}\label{new-eq-3}
 (x)_k=(-1)^k(x+j_1)\cdots (x+j_i)\frac{\Gamma_p(x+k)}{\Gamma_p(x)}.
\end{align}
We note that \eqref{new-eq-2} and \eqref{new-eq-3} follow from part (2) of Proposition \ref{prop-1} for every $x\in \mathbb{Z}_p$.
\par
We now state a product formula for the $p$-adic Gamma function from \cite[p. 5]{mccarthy2}.  Let $\mathbb{F}_p$ denote the finite field with $p$ elements.
Let $\omega: \mathbb{F}_p^{\times} \rightarrow \mathbb{Z}_p^{\times}$ be the Teichm\"{u}ller character.
For $a\in\mathbb{F}_p^\times$, the value $\omega(a)$ is just the $(p-1)$-th root of unity in $\mathbb{Z}_p$
such that $\omega(a)\equiv a \pmod{p}$. If $m\in\mathbb{Z}^+$, $p\nmid m$ and $x =\frac{r}{p-1}$ with $0 \leq r \leq p-1$, then
\begin{align}\label{prod-formula}
\prod_{h=0}^{m-1}\Gamma_p\left(\frac{x+h}{m}\right)=\omega \left( m^{(1-x)(1-p)}\right)
\Gamma_p(x)
\prod_{h=1}^{m-1}\Gamma_p\left(\frac{h}{m}\right).
\end{align}
Finally, we recall some hypergeometric formulas from \cite{andrews, bailey, gasper, slater}.
We first state the Whipple's well-posed ${_{7}}F_{6}$ evaluation formula.
\begin{theorem}\cite[Theorem 3.4.5]{andrews}\label{thm-4}
We have
\begin{align}
&{_{7}}F_{6}\left[\begin{array}{ccccccc}
                 a, & \frac{a}{2}+1, & b, & c, & d, & e, & f \\
                 ~ & \frac{a}{2}, & 1+a-b, & 1+a-c, & 1+a-d, & 1+a-e, & 1+a-f
               \end{array};1
\right]\notag\\
&=\frac{\Gamma(1+a-d)\Gamma(1+a-e)\Gamma(1+a-f)\Gamma(1+a-d-e-f)}{\Gamma(1+a)\Gamma(1+a-e-f)\Gamma(1+a-d-e)\Gamma(1+a-d-f)}\notag\\
&\hspace{2cm}\times {_4}F_3\left[\begin{array}{cccc}
                                                                      1+a-b-c, & d, & e, & f \\
                                                                     ~ & 1+a-b, & 1+a-c, & d+e+f-a
                                                                    \end{array}; 1
\right],\notag
\end{align}
provided the left side converges and the right side terminates.
\end{theorem}
We next state the Dougall's formula.
\begin{theorem}\cite[Theorem 3.5.1]{andrews}\label{thm-2}
If $f$ is a negative integer and $1+2a=b+c+d+e+f$ then
\begin{align}
&{_{7}}F_{6}\left[\begin{array}{ccccccc}
                 a, & \frac{a}{2}+1, & b, & c, & d, & e, & f \\
                 ~ & \frac{a}{2}, & 1+a-b, & 1+a-c, & 1+a-d, & 1+a-e, & 1+a-f
               \end{array};1
\right]\notag\\
&=\frac{(1+a)_{-f}(1+a-b-c)_{-f}(1+a-b-d)_{-f}(1+a-c-d)_{-f}}{(1+a-b)_{-f}(1+a-c)_{-f}(1+a-d)_{-f}(1+a-b-c-d)_{-f}}.\notag
\end{align}
\end{theorem}
We also need the following identity to prove our main results.
\begin{theorem}\cite[p. 56]{slater}\label{thm-6}
If $m$ is a positive integer then
\begin{align}
&{_5}F_{4}\left[\begin{array}{ccccc}
                  a, & 1+\frac{a}{2}, & b, & c, & -m \\
                  ~ & \frac{a}{2}, & 1+a-b, & 1+a-c, & 1+a+m
                \end{array};1
\right]
=\frac{(1+a)_{m}(1+a-b-c)_{m}}{(1+a-b)_{m}(1+a-c)_{m}}.\notag
\end{align}
\end{theorem}
We now state the following particular case of the Karlsson-Minton formula.
\begin{theorem}\cite[Eqn. (1.9.3)]{gasper}\label{thm-7}
For non-negative integers $m_1,m_2,\ldots m_n$ we have
\begin{align}
&{_{n+1}}F_n\left[\begin{array}{cccc}
  -(m_1+m_2+\cdots +m_n), & b_1+m_1, & \ldots, & b_n+m_n  \\
  ~ & b_1, & \ldots, & b_n
\end{array};1\right]\notag\\
&=(-1)^{m_1+m_2+\cdots+m_n}\frac{(m_1+m_2+\cdots+m_n)!}{(b_1)_{m_1}\cdots (b_n)_{m_n}}.\notag
\end{align}
\end{theorem}
\section{Proof of the results}
Various supercongruences have been proved by many mathematicians using a variety of methods.
Recently, L. Long and R. Ramakrishna \cite{long}, A. Deines et al. \cite{DFLST} and H. Swisher \cite{swisher} prove several supercongruences
using a technique which relies on the relations between the classical and $p$-adic Gamma functions.
We use a similar technique. In the proof, we consider certain
classical hypergeometric evaluation identities and then use the relation between the classical and $p$-adic Gamma function.
\begin{proof}[Proof of Theorem \ref{MT-5}]
Let $m_1=m_2=\cdots=m_{n-1}=\frac{p-1}{n(n-1)}$, $b_1=\frac{1}{n-1}+\frac{yp}{n-1}$ and $b_2=\cdots=b_{n-1}=\frac{1}{n-1}$, where 
$y\in \mathbb{Z}_p$. Substituting
these values in Theorem \ref{thm-7} we have
\begin{align}\label{eq-32}
&{_nF_{n-1}}\left[\begin{array}{ccccc}
                    \frac{1-p}{n}, & \frac{1}{n}+\frac{p}{n-1}(y+\frac{1}{n}), & \frac{1}{n}+\frac{p}{n(n-1)}, & \ldots,
                     & \frac{1}{n}+\frac{p}{n(n-1)} \vspace{.1cm}\\
                    ~ & \frac{1}{n-1}+\frac{yp}{n-1}, & \frac{1}{n-1}, & \ldots,  & \frac{1}{n-1}
                  \end{array};1
\right]_\frac{p-1}{n}\notag\\
&=\dfrac{(-1)^{\frac{p-1}{n}}(\frac{p-1}{n})!}{(\frac{1}{n-1}+\frac{yp}{n-1})_{\frac{p-1}{n(n-1)}}
(\frac{1}{n-1})_{\frac{p-1}{n(n-1)}}\cdots (\frac{1}{n-1})_{\frac{p-1}{n(n-1)}}}.
\end{align}
By definition we have
\begin{align}\label{eq-30}
&{_{n}}F_{n-1}\left[\begin{array}{ccccc}
                    \frac{1-p}{n}, & \frac{1}{n}+\frac{p}{n-1}(y+\frac{1}{n}), & \frac{1}{n}+\frac{p}{n(n-1)}, & \ldots,
                     & \frac{1}{n}+\frac{p}{n(n-1)} \vspace{.1cm}\\
                    ~ & \frac{1}{n-1}+\frac{yp}{n-1}, & \frac{1}{n-1}, & \ldots,  & \frac{1}{n-1}
                  \end{array};1
\right]_\frac{p-1}{n}\notag\\
&=\sum_{k=0}^{\frac{p-1}{n}}\frac{(\frac{1-p}{n})_k~(\frac{1}{n}+\frac{p}{n-1}(y+\frac{1}{n}))_k~
(\frac{1}{n}+\frac{p}{n(n-1)})_k^{n-2}}{(\frac{1}{n-1}+\frac{yp}{n-1})_k~(\frac{1}{n-1})_k^{n-2}}\frac{1}{k!}.
\end{align}
Now,
\begin{align}
\left(\frac{1-p}{n}\right)_k&=\left(\frac{1}{n}-\frac{p}{n}\right)\left(\frac{1}{n}+1-\frac{p}{n}\right)\cdots\left(\frac{1}{n}+k-1-\frac{p}{n}\right)\notag\\
&=\frac{1}{n}\left(\frac{1}{n}+1\right)\left(\frac{1}{n}+2\right)\cdots\left(\frac{1}{n}+k-1\right)\prod_{i=0}^{k-1}
\left(1-\frac{p}{n(\frac{1}{n}+i)}\right)\notag\\
&\equiv\left(\frac{1}{n}\right)_k\left(1-\frac{p}{n}\sum_{i=0}^{k-1}\frac{1}{\frac{1}{n}+i}\right)\pmod {p^2}\notag\\
&\equiv\left(\frac{1}{n}\right)_k\left(1-\frac{p}{n}A_k\right)\pmod {p^2},\notag
\end{align}
where $A_k=\displaystyle\sum_{i=0}^{k-1}\frac{1}{\frac{1}{n}+i}\in\mathbb{Z}_p$.
Similarly, we have
\begin{align}
\left(\frac{1}{n}+\frac{p}{n-1}\left(y+\frac{1}{n}\right)\right)_k\equiv\left(\frac{1}{n}\right)_k\left(1+\frac{p(\frac{1}{n}+y)}{n-1}
A_k\right)\pmod {p^2},\notag
\end{align}
\begin{align}
\left(\frac{1}{n}+\frac{p}{n(n-1)}\right)_k^{n-2}&\equiv\left(\frac{1}{n}\right)_k^{n-2}\left(1+\frac{p}{n(n-1)}
A_k\right)^{n-2}\pmod {p^2}\notag\\
&\equiv\left(\frac{1}{n}\right)_k^{n-2}\left(1+\frac{p(n-2)}{n(n-1)}
A_k\right)\pmod {p^2}\notag
\end{align}
and
\begin{align}
\left(\frac{1}{n-1}+\frac{yp}{n-1}\right)_k^{-1}&\equiv\left(\frac{1}{n-1}\right)_k^{-1}
\left(1+\frac{yp}{n-1}\sum_{i=0}^{k-1}\frac{1}{\frac{1}{n-1}+i}\right)^{-1}\pmod {p^2}\notag\\
&\equiv\left(\frac{1}{n-1}\right)_k^{-1}
\left(1-\frac{yp}{n-1}B_k\right)\pmod {p^2},\notag
\end{align}
where $B_k=\displaystyle\sum_{i=0}^{k-1}\frac{1}{\frac{1}{n-1}+i}\in\mathbb{Z}_p$.
Substituting all these identities on the right hand side of \eqref{eq-30} we deduce that
\begin{align}\label{eq-31}
&{_{n}}F_{n-1}\left[\begin{array}{ccccc}
                    \frac{1-p}{n}, & \frac{1}{n}+\frac{p}{n-1}(y+\frac{1}{n}), & \frac{1}{n}+\frac{p}{n(n-1)}, & \ldots,
                     & \frac{1}{n}+\frac{p}{n(n-1)} \vspace{.1cm}\\
                    ~ & \frac{1}{n-1}+\frac{yp}{n-1}, & \frac{1}{n-1}, & \ldots,  & \frac{1}{n-1}
                  \end{array};1
\right]_\frac{p-1}{n}\notag\\
&\equiv\sum_{k=0}^{\frac{p-1}{n}}\frac{(\frac{1}{n})_k^{n}}{(\frac{1}{n-1})_k^{n-1}k!}\left\{1+\frac{yp}{n-1}(A_k-B_k)\right\}
\pmod{p^2}\notag\\
&\equiv{_n}F_{n-1}\left[\begin{array}{cccc}
                          \frac{1}{n}, & \frac{1}{n}, & \ldots, & \frac{1}{n} \vspace{.1cm}\\
                          ~ & \frac{1}{n-1}, & \ldots, & \frac{1}{n-1}
                        \end{array};1
\right]_{\frac{p-1}{n}}+\frac{yp}{n-1}C\pmod{p^2}
\end{align}
for some $C\in \mathbb{Z}_p$.
We know that $\left(\frac{p-1}{n}\right)!=(1)_{\frac{p-1}{n}}$. Now, using \eqref{new-eq-2} and the fact that $\Gamma_p(1)=-1$ 
we obtain $$\left(\frac{p-1}{n}\right)!=(1)_{\frac{p-1}{n}}=(-1)^{1+\frac{p-1}{n}}\Gamma_p\left(1+\frac{p-1}{n}\right)=
(-1)^{1+\frac{p-1}{n}}\Gamma_p\left(\frac{n-1}{n}+\frac{p}{n}\right).$$ Again, \eqref{new-eq-2}
yields
\begin{align}
\left(\frac{1}{n-1}+\frac{yp}{n-1}\right)_{\frac{p-1}{n(n-1)}}&=(-1)^{\frac{p-1}{n(n-1)}}
\frac{\Gamma_p(\frac{1}{n-1}+\frac{yp}{n-1}+\frac{p-1}{n(n-1)})}{\Gamma_p(\frac{1}{n-1}+\frac{yp}{n-1})}\notag\\
&=(-1)^{\frac{p-1}{n(n-1)}}\frac{\Gamma_p(\frac{1}{n}+\frac{p(1+ny)}{n(n-1)})}{\Gamma_p(\frac{1}{n-1}+\frac{yp}{n-1})},\notag
\end{align}
and
\begin{align}
\left(\frac{1}{n-1}\right)_{\frac{p-1}{n(n-1)}}=(-1)^{\frac{p-1}{n(n-1)}}\frac{\Gamma_p(\frac{1}{n-1}+\frac{p-1}{n(n-1)})}
{\Gamma_p(\frac{1}{n-1})}
=(-1)^{\frac{p-1}{n(n-1)}}\frac{\Gamma_p(\frac{1}{n}+\frac{p}{n(n-1)})}
{\Gamma_p(\frac{1}{n-1})}.\notag
\end{align}
Putting these values in \eqref{eq-32}, we find that the right hand side reduces to
\begin{align}\label{rv-eq-1}
\frac{(-1)^{1+\frac{p-1}{n}}~\Gamma_p(\frac{n-1}{n}+\frac{p}{n})~\Gamma_p(\frac{1}{n-1}+\frac{yp}{n-1})~\Gamma_p(\frac{1}{n-1})^{n-2}}
{\Gamma_p(\frac{1}{n}+\frac{p(1+ny)}{n(n-1)})~\Gamma_p(\frac{1}{n}+\frac{p}{n(n-1)})^{n-2}}.
\end{align}
Applying Theorem \ref{thm-long-14} for $p\geq 5$ we have
\begin{align}
\Gamma_p\left(\frac{n-1}{n}+\frac{p}{n}\right)\equiv\Gamma_p\left(\frac{n-1}{n}\right)\left(1+G_1\left(\frac{n-1}{n}\right)\frac{p}{n}\right)\pmod{p^2}.\notag
\end{align}
From \cite[Proposition 3.2 (3)]{mccarthy} we have $G_1\left(\frac{1}{n}\right)=G_1\left(\frac{n-1}{n}\right)$. Using this value we obtain
\begin{align}\label{rv-eq-2}
\Gamma_p\left(\frac{n-1}{n}+\frac{p}{n}\right)\equiv\Gamma_p\left(\frac{n-1}{n}\right)\left(1+G_1\left(\frac{1}{n}\right)\frac{p}{n}\right)\pmod{p^2}.
\end{align}
Again, applying Theorem \ref{thm-long-14} for $p\geq 5$ we have
\begin{align}\label{rv-eq-3}
\Gamma_p\left(\frac{1}{n-1}+\frac{yp}{n-1}\right)\equiv\Gamma_p\left(\frac{1}{n-1}\right)\left(1+G_1\left(\frac{1}{n-1}\right)\frac{yp}{n-1}\right)\pmod{p^2},
\end{align}
\begin{align}\label{rv-eq-4}
\Gamma_p\left(\frac{1}{n}+\frac{p(1+ny)}{n(n-1)}\right)^{-1}&\equiv\Gamma_p\left(\frac{1}{n}\right)^{-1}
\left(1+G_1\left(\frac{1}{n}\right)\frac{p(1+ny)}{n(n-1)}\right)^{-1}\pmod{p^2}\notag\\
&\equiv\Gamma_p\left(\frac{1}{n}\right)^{-1}
\left(1-G_1\left(\frac{1}{n}\right)\frac{p(1+ny)}{n(n-1)}\right)\pmod{p^2},
\end{align}
and
\begin{align}\label{rv-eq-5}
\Gamma_p\left(\frac{1}{n}+\frac{p}{n(n-1)}\right)\equiv\Gamma_p\left(\frac{1}{n}\right)\left(1+G_1\left(\frac{1}{n}\right)\frac{p}{n(n-1)}\right)\pmod{p^2}.
\end{align}

If we use \eqref{rv-eq-1}, \eqref{rv-eq-2}, \eqref{rv-eq-3}, \eqref{rv-eq-4} and \eqref{rv-eq-5}, then modulo $p^2$ the right hand side of \eqref{eq-32} becomes
\begin{align}\label{eq-33}
\dfrac{(-1)^{1+\frac{p-1}{n}}\Gamma_p(\frac{n-1}{n})\Gamma_p(\frac{1}{n-1})^{n-1}}{\Gamma_p(\frac{1}{n})^{n-1}}
\left[1+\left(G_1\left(\frac{1}{n-1}\right)-G_1\left(\frac{1}{n}\right)\right)\frac{yp}{n-1}\right].
\end{align}
Since $p\equiv 1 \pmod{n(n-1)}$, from the $p$-adic expansion of $-1$ we have $$-\frac{1}{n}=\frac{p-1}{n}+\frac{p-1}{n}p+\cdots. $$
Hence, $a_0(\frac{1}{n})=p-\frac{p-1}{n}$ and part (3) of Proposition \ref{prop-1} yields
\begin{align}\label{eqn-new-1}
\displaystyle \frac{1}{\Gamma_p(\frac{1}{n})^{n-1}}&=(-1)^{(n-1)p}(-1)^{\frac{(n-1)(p-1)}{n}}\displaystyle ~\Gamma_p\left(\frac{n-1}{n}\right)^{n-1}\notag\\
&=(-1)^{(n-1)+\frac{p-1}{n}}~\Gamma_p\left(\frac{n-1}{n}\right)^{n-1}.
\end{align}
Now, from \eqref{eq-32}, \eqref{eq-31}, \eqref{eq-33} and \eqref{eqn-new-1} we deduce that, modulo $p^2$
\begin{align}\label{eq-34}
&{_n}F_{n-1}\left[\begin{array}{cccc}
                          \frac{1}{n}, & \frac{1}{n}, & \ldots, & \frac{1}{n} \vspace{.1cm}\\
                          ~ & \frac{1}{n-1}, & \ldots, & \frac{1}{n-1}
                        \end{array};1
\right]_{\frac{p-1}{n}}+\frac{yp}{n-1}C\notag \\
&\equiv(-1)^n\Gamma_p\left(\frac{1}{n-1}\right)^{n-1}\Gamma_p\left(\frac{n-1}{n}\right)^n\left[1+\left(G_1\left(\frac{1}{n-1}\right)
-G_1\left(\frac{1}{n}\right)\right)\frac{yp}{n-1}\right],
\end{align}
which is true for all $y\in \mathbb{Z}_p$. Putting $y=0$ in \eqref{eq-34} we obtain that, modulo $p^2$ 
\begin{align}\label{rv-eq-6}
{_n}F_{n-1}\left[\begin{array}{cccc}
                          \frac{1}{n}, & \frac{1}{n}, & \ldots, & \frac{1}{n} \vspace{.1cm}\\
                          ~ & \frac{1}{n-1}, & \ldots, & \frac{1}{n-1}
                        \end{array};1
\right]_{\frac{p-1}{n}}\equiv(-1)^n\Gamma_p\left(\frac{1}{n-1}\right)^{n-1}\Gamma_p\left(\frac{n-1}{n}\right)^n.
\end{align}
Since for $\frac{p-1}{n}<k\leq\frac{p-1}{n-1}$, we know that $\left(\frac{1}{n}\right)_k$ contains a factor $\frac{p}{n}$ and 
$\left(\frac{1}{n-1}\right)_k$ and $k!$ have no factor containing $p$, hence modulo $p^2$ we have
\begin{align}\label{rv-eq-7}
{_n}F_{n-1}\left[\begin{array}{cccc}
                          \frac{1}{n}, & \frac{1}{n}, & \ldots, & \frac{1}{n} \vspace{.1cm}\\
                          ~ & \frac{1}{n-1}, & \ldots, & \frac{1}{n-1}
                        \end{array};1
\right]_{\frac{p-1}{n}}\equiv{_n}F_{n-1}\left[\begin{array}{cccc}
                          \frac{1}{n}, & \frac{1}{n}, & \ldots, & \frac{1}{n} \vspace{.1cm}\\
                          ~ & \frac{1}{n-1}, & \ldots, & \frac{1}{n-1}
                        \end{array};1
\right]_{\frac{p-1}{n-1}}.
\end{align}
Now, combining \eqref{rv-eq-6} and \eqref{rv-eq-7} we obtain the result. We have verified the case $p=3$ by hand.
This completes the proof of the theorem.
\end{proof}
\begin{proof}[Proof of Theorem \ref{MT-4}]
Let $\zeta$ be a primitive cubic root of unity.
Putting $a=\frac{1}{n}$, $b=\frac{1-\zeta p}{n}$, $c=\frac{1-\zeta^2p}{n}$ and $m=\frac{p-1}{n}$ in Theorem \ref{thm-6} we have
\begin{align}\label{eq-25}
&{_5}F_4\left[\begin{array}{ccccc}
               \frac{1}{n}, & 1+\frac{1}{2n}, & \frac{1-\zeta p}{n}, & \frac{1-\zeta^2p}{n}, & \frac{1-p}{n} \\
               ~ & \frac{1}{2n}, & 1+\frac{\zeta p}{n}, & 1+\frac{\zeta^2p}{n}, & 1+\frac{p}{n}
             \end{array};1
\right]_{\frac{p-1}{n}}\notag\\
&=\frac{(1+\frac{1}{n})_{\frac{p-1}{n}}(\frac{n-1}{n}+\frac{(\zeta+\zeta^2)p}{n})_{\frac{p-1}{n}}}
{(1+\frac{\zeta p}{n})_{\frac{p-1}{n}}(1+\frac{\zeta^2p}{n})_{\frac{p-1}{n}}}.
\end{align}
We first show that
\begin{align}\label{eq-770}
&{_5}F_4\left[\begin{array}{ccccc}
               \frac{1}{n}, & 1+\frac{1}{2n}, & \frac{1-\zeta p}{n}, & \frac{1-\zeta^2p}{n}, & \frac{1-p}{n} \\
               ~ & \frac{1}{2n}, & 1+\frac{\zeta p}{n}, & 1+\frac{\zeta^2p}{n}, & 1+\frac{p}{n}
             \end{array};1
\right]_{\frac{p-1}{n}}\notag\\
&\equiv{_5}F_4\left[\begin{array}{ccccc}
               \frac{1}{n}, & 1+\frac{1}{2n}, & \frac{1}{n}, & \frac{1}{n}, & \frac{1}{n} \\
               ~ & \frac{1}{2n}, & 1, & 1, & 1
             \end{array};1
\right]_{\frac{p-1}{n}}\notag\\&\equiv
{_5}F_4\left[\begin{array}{ccccc}
               \frac{1}{n}, & 1+\frac{1}{2n}, & \frac{1}{n}, & \frac{1}{n}, & \frac{1}{n} \\
               ~ & \frac{1}{2n}, & 1, & 1, & 1
             \end{array};1
\right]_{p-1}\pmod{p^4}.
\end{align}
The last congruence follows from the fact that the rising factorial $(\frac{1}{n})_k$ contains a multiple of $p$
for $\frac{p-1}{n}<k\leq p-1$.
\par Now,
\begin{align*}
 \displaystyle \frac{(\frac{1}{n}-\frac{x}{n})_k}{(1+\frac{x}{n})_k}=\prod_{j=0}^{k-1}\frac{(\frac{1}{n}+j-\frac{x}{n})}{(1+\frac{x}{n}+j)}=
 \prod_{j=0}^{k-1}\frac{(\frac{1}{n}+j)(1-\frac{x/n}{1/n+j})}{(1+j)(1+\frac{x/n}{1+j})}.
\end{align*}
The terms $\frac{1}{n}+j$ and $1+j$ do not contain a multiple of $p$ for $0\leq k\leq \frac{p-1}{n}$. Therefore, we can find $a_{k,1},a_{k,2}, \ldots\in\mathbb{Z}_p$ such that
\begin{align}\label{rev-new-eq-1}
\dfrac{(\frac{1}{n}-\frac{x}{n})_k}{(1+\frac{x}{n})_k}&=\dfrac{(\frac{1}{n})_k}{(1)_k}[1+\sum_{i\geq1}a_{k,i}x^i]
\end{align}
for $0\leq k\leq\frac{p-1}{n}$. Thus,
\begin{align}
&{_5}F_4\left[\begin{array}{ccccc}
               \frac{1}{n}, & 1+\frac{1}{2n}, & \frac{1-x}{n}, & \frac{1-y}{n}, & \frac{1-p}{n} \\
               ~ & \frac{1}{2n}, & 1+\frac{x}{n}, & 1+\frac{y}{n}, & 1+\frac{p}{n}
             \end{array};1
\right]_{\frac{p-1}{n}}\notag\\&\equiv
{_5}F_4\left[\begin{array}{ccccc}
               \frac{1}{n}, & 1+\frac{1}{2n}, & \frac{1-x}{n}, & \frac{1-y}{n}, & \frac{1}{n} \\
               ~ & \frac{1}{2n}, & 1+\frac{x}{n}, & 1+\frac{y}{n}, & 1
             \end{array};1
\right]_{\frac{p-1}{n}}\pmod{p}\notag\\
&=\sum_{k=0}^{\frac{p-1}{n}}\dfrac{(1+\frac{1}{2n})_k(\frac{1}{n})_k^4}
{(\frac{1}{2n})_kk!^4}\left[1+\sum_{i\geq1}a_{k,i}x^i\right]\left[1+\sum_{i\geq1}a_{k,i}y^i\right]\in\mathbb{Z}_p[[x,y]].\notag
\end{align}
Again, Theorem \ref{thm-6} implies that
\begin{align}
&{_5}F_4\left[\begin{array}{ccccc}
               \frac{1}{n}, & 1+\frac{1}{2n}, & \frac{1-x}{n}, & \frac{1-y}{n}, & \frac{1-p}{n} \\
               ~ & \frac{1}{2n}, & 1+\frac{x}{n}, & 1+\frac{y}{n}, & 1+\frac{p}{n}
             \end{array};1
\right]_{\frac{p-1}{n}}\notag\\
&=\frac{(1+\frac{1}{n})_{\frac{p-1}{n}}(\frac{n-1}{n}+\frac{(x+y)}{n})_{\frac{p-1}{n}}}
{(1+\frac{x}{n})_{\frac{p-1}{n}}(1+\frac{y}{n})_{\frac{p-1}{n}}}.\notag
\end{align}
We note that the rising factorial $(1+\frac{1}{n})_{\frac{p-1}{n}}$ contains a factor which is a multiple of $p$, namely
$(1+\frac{1}{n}+\frac{p-1}{n}-1)=\frac{p}{n}$; and the other factors do not contain any multiple of $p$. Therefore, we have
\begin{align}\label{eq-721}
&{_5}F_4\left[\begin{array}{ccccc}
               \frac{1}{n}, & 1+\frac{1}{2n}, & \frac{1-x}{n}, & \frac{1-y}{n}, & \frac{1-p}{n} \\
               ~ & \frac{1}{2n}, & 1+\frac{x}{n}, & 1+\frac{y}{n}, & 1+\frac{p}{n}
             \end{array};1
\right]_{\frac{p-1}{n}}\in p\mathbb{Z}_p[[x,y]].
\end{align}
If we put $x=\zeta p$ and $y=\zeta^2 p$, then using Lemma \ref{lemma-swisher} we obtain
\begin{align}
&{_5}F_4\left[\begin{array}{ccccc}
               \frac{1}{n}, & 1+\frac{1}{2n}, & \frac{1-\zeta p}{n}, & \frac{1-\zeta^2 p}{n}, & \frac{1-p}{n} \\
               ~ & \frac{1}{2n}, & 1+\frac{\zeta p}{n}, & 1+\frac{\zeta^2 p}{n}, & 1+\frac{p}{n}
             \end{array};1
\right]_{\frac{p-1}{n}}\notag\\&
=A_0+\sum_{i\geq1}A_{3i}p^{3i},\notag
\end{align}
where $\displaystyle A_0=\sum_{k=0}^{\frac{p-1}{n}}\dfrac{(1+\frac{1}{2n})_k(\frac{1}{n})_k^4}
{(\frac{1}{2n})_kk!^4}$ and $\displaystyle A_{3i}=\sum_{k=0}^{\frac{p-1}{n}}\dfrac{(1+\frac{1}{2n})_k(\frac{1}{n})_k^4}
{(\frac{1}{2n})_kk!^4}H(a_{k,i})$ and $H(a_{k,i})$ is an integral polynomial in the $a_{k,i}$ where the
second subscripts in each monomial add to $3$. Now, \eqref{eq-721} implies that $A_0$ and $A_{3i}\in p\mathbb{Z}_p$.
Thus,
\begin{align}\label{eq-28}
&{_5}F_4\left[\begin{array}{ccccc}
               \frac{1}{n}, & 1+\frac{1}{2n}, & \frac{1-\zeta p}{n}, & \frac{1-\zeta^2p}{n}, & \frac{1-p}{n} \\
               ~ & \frac{1}{2n}, & 1+\frac{\zeta p}{n}, & 1+\frac{\zeta^2p}{n}, & 1+\frac{p}{n}
             \end{array};1
\right]_{\frac{p-1}{n}}\notag\\
&\equiv\sum_{k=0}^{\frac{p-1}{n}}\dfrac{(1+\frac{1}{2n})_k(\frac{1}{n})_k^4}
{(\frac{1}{2n})_kk!^4}\notag\\
&\equiv{_5}F_4\left[\begin{array}{ccccc}
               \frac{1}{n}, & 1+\frac{1}{2n}, & \frac{1}{n}, & \frac{1}{n}, & \frac{1}{n} \\
               ~ & \frac{1}{2n}, & 1, & 1, & 1
             \end{array};1
\right]_{\frac{p-1}{n}}\pmod {p^4},
\end{align}
which proves \eqref{eq-770}.
\par
Since the rising factorial $(1+\frac{1}{n})_{\frac{p-1}{n}}$ contains a factor which is a multiple of $p$, namely
$(1+\frac{1}{n}+\frac{p-1}{n}-1)=\frac{p}{n}$, using \eqref{new-eq-3} and then part (2) of Proposition \ref{prop-1}, we deduce that
\begin{align*}
 \displaystyle \left(1+\frac{1}{n}\right)_{\frac{p-1}{n}}=(-1)^{\frac{p-1}{n}}\frac{p}{n}\frac{\Gamma_p(1+\frac{p}{n})}{\Gamma_p(1+\frac{1}{n})}
 =-(-1)^{\frac{p-1}{n}}p\frac{\Gamma_p(1+\frac{p}{n})}{\Gamma_p(\frac{1}{n})}.
\end{align*}
Also, applying \eqref{new-eq-2} for the remaining three rising factorials on the right hand side of \eqref{eq-25},  we deduce that
\begin{align}\label{eq-29}
&\dfrac{(1+\frac{1}{n})_{\frac{p-1}{n}}(\frac{n-1}{n}+\frac{(\zeta+\zeta^2)p}{n})_{\frac{p-1}{n}}}
{(1+\frac{\zeta p}{n})_{\frac{p-1}{n}}(1+\frac{\zeta^2p}{n})_{\frac{p-1}{n}}}\notag\\
&=-p\cdot \dfrac{\Gamma_p(1+\frac{p}{n})\Gamma_p(1+\frac{p\zeta}{n})\Gamma_p(1+\frac{p\zeta^2}{n})\Gamma_p(\frac{n-2}{n})}
{\Gamma_p(\frac{1}{n})\Gamma_p(\frac{n-1}{n}-\frac{p}{n})\Gamma_p(\frac{n-1}{n}-\frac{p\zeta}{n})
\Gamma_p(\frac{n-1}{n}-\frac{p\zeta^2}{n})}.
\end{align}
Now, using Theorem \ref{thm-long-14}, we find that
\begin{align}\label{eqn-new-3}
\Gamma_p\left(1+\frac{p}{n}\right)\Gamma_p\left(1+\frac{\zeta p}{n}\right)\Gamma_p\left(1+\frac{\zeta^2p}{n}\right)
&=\Gamma_p(1)^3(1+O(p^3))=-1+O(p^3).
\end{align}
Similarly, using the last part of Proposition \ref{prop-2} and then Theorem \ref{thm-long-14}, we find that
\begin{align}\label{eqn-new-4}
&\frac{1}{\Gamma_p(1-\frac{1}{n}-\frac{p}{n})\Gamma_p(1-\frac{1}{n}-\frac{\zeta p}{n})\Gamma_p(1-\frac{1}{n}-\frac{\zeta^2p}{n})}\notag\\
&=(-1)^{a_0(\frac{1}{n})}\Gamma_p\left(\frac{1}{n}\right)^3\left(1+O(p^3)\right)\notag\\
&=(-1)^{1+\frac{p-1}{n}}\Gamma_p\left(\frac{1}{n}\right)^3\left(1+O(p^3)\right).
\end{align}
From \eqref{eq-29}, \eqref{eqn-new-3} and \eqref{eqn-new-4}, we have, for $p\geq 5$
\begin{align}\label{new-eq-29}
&\dfrac{(1+\frac{1}{n})_{\frac{p-1}{n}}(\frac{n-1}{n}+\frac{(\zeta+\zeta^2)p}{n})_{\frac{p-1}{n}}}
{(1+\frac{\zeta p}{n})_{\frac{p-1}{n}}(1+\frac{\zeta^2p}{n})_{\frac{p-1}{n}}}\notag\\
&\equiv (-1)^{1+\frac{p-1}{n}} p\Gamma_p\left(\frac{1}{n}\right)^2\Gamma_p\left(\frac{n-2}{n}\right)\pmod {p^4}.
\end{align}
Finally, from \eqref{eq-25}, \eqref{eq-770}, \eqref{eq-28} and \eqref{new-eq-29} we obtain our result.
\end{proof}
\begin{proof}[Proof of Theorem \ref{MT-6}]
Let $a=\frac{1}{n}$, $b=\frac{2n-3}{2n}$,
$c=\frac{1}{2n}$, $d=\frac{1-\zeta p}{n}$, $e=\frac{1-\zeta^2p}{n}$ and
$f=\frac{1-p}{n}$, where $\zeta$ is a primitive cubic root of unity.
Plugging these values in Theorem \ref{thm-2} we have
\begin{align}\label{eq-200}
&_{7}F_6\left[\begin{array}{ccccccc}
                 \frac{1}{n}, & \frac{2n+1}{2n}, & \frac{2n-3}{2n}, & \frac{1}{2n}, & \frac{1-\zeta p}{n}, & \frac{1-\zeta^2p}{n},
                 & \frac{1-p}{n} \\
                 ~ & \frac{1}{2n}, & \frac{5}{2n}, & \frac{2n+1}{2n}, & 1+\frac{\zeta p}{n}, & 1+\frac{\zeta^2p}{n}, & 1+\frac{p}{n}
               \end{array};1
\right]_{\frac{p-1}{n}}\notag\\
&={_5}F_4\left[\begin{array}{ccccc}
                \frac{1}{n}, & \frac{2n-3}{2n}, & \frac{1-\zeta p}{n}, & \frac{1-\zeta^2p}{n}, & \frac{1-p}{n} \\
                ~ & \frac{5}{2n}, & 1+\frac{\zeta p}{n}, & 1+\frac{\zeta^2p}{n}, & 1+\frac{p}{n}
              \end{array};1
\right]_{\frac{p-1}{n}}\notag\\
&=\dfrac{(\frac{n+1}{n})_{\frac{p-1}{n}}(\frac{2}{n})_{\frac{p-1}{n}}(\frac{3}{2n}+\frac{\zeta p}{n})_{\frac{p-1}{n}}(\frac{2n-1}{2n}
+\frac{\zeta p}{n})_{\frac{p-1}{n}}}
{(\frac{5}{2n})_{\frac{p-1}{n}}(\frac{2n+1}{2n})_{\frac{p-1}{n}}(1+\frac{\zeta p}{n})_{\frac{p-1}{n}}(\frac{1}{n}+\frac{\zeta p}{n})_{\frac{p-1}{n}}}.
\end{align}
We now prove that
\begin{align}\label{eq-400}
&_{5}F_4\left[\begin{array}{ccccc}
                \frac{1}{n}, & \frac{2n-3}{2n}, & \frac{1-\zeta p}{n}, & \frac{1-\zeta^2p}{n}, & \frac{1-p}{n} \\
                ~ & \frac{5}{2n}, & 1+\frac{\zeta p}{n}, & 1+\frac{\zeta^2p}{n}, & 1+\frac{p}{n}
              \end{array};1
\right]_{\frac{p-1}{n}}\notag\\
&\equiv
{_5}F_{4}\left[\begin{array}{ccccc}
                 \frac{1}{n}, & \frac{1}{n}, & \frac{1}{n}, & \frac{1}{n}, & \frac{2n-3}{2n} \vspace{.08cm}\\
                 ~ & \frac{5}{2n}, & 1, & 1, & 1
               \end{array};1
\right]_{\frac{p-1}{n}}\notag\\&
\equiv{_5}F_{4}\left[\begin{array}{ccccc}
                 \frac{1}{n}, & \frac{1}{n}, & \frac{1}{n}, & \frac{1}{n}, & \frac{2n-3}{2n} \vspace{.08cm}\\
                 ~ & \frac{5}{2n}, & 1, & 1, & 1
               \end{array};1
\right]_{p-1}\pmod{p^4}.
\end{align}
The last congruence follows from the fact that the rising factorial $(\frac{1}{n})_k$ contains a multiple of $p$ for $\frac{p-1}{n}<k\leq p-1$. Further, if $(\frac{5}{2n})_k$
contains a multiple of $p$, then the term $(\frac{2n-3}{2n})_k$ also contains a multiple of $p$. For example, let 
$\frac{5}{2n}+\ell=\frac{5+2n\ell}{2n}=\frac{\alpha p}{2n}$ for some $\alpha$, where $\ell=0, 1, \ldots, k-1$. 
Thus $(\frac{5}{2n})_k$ contains a multiple of $p$ for $k=\frac{\alpha p-5}{2n}+1, \frac{\alpha p-5}{2n}+2,\ldots, p-1$.
This implies that $\alpha p\equiv5\pmod{2n}$, and hence $\alpha p\equiv5\pmod{n}$. 
Also, we have $p\equiv1\pmod{n}$, which gives $\alpha\equiv5\pmod{n}$. Since $n$ is even, $\alpha$ must be odd. 
    Now, consider $k_1=\frac{\alpha p-5}{2n}-\frac{4(p-1)}{n}=\frac{(\alpha-8)p+3}{2n}$. 
    Then $0<k_1<k\leq p-1$ for $\alpha>8$ and $(\frac{2n-3}{2n})_{k_1}$ contains a multiple of $p$. 
    The other remaining odd values of $\alpha$ are $1, 3, 5, 7$.
    Since $n\neq2,4$, we easily find that $\alpha \neq 1, 3, 7$. Now let $\alpha=5$. If $\gcd(5, n)=1$ then from $5p\equiv5\pmod{2n}$ we have $p\equiv 1\pmod{2n}$, which 
    is not possible as $p\not \equiv 1\pmod{2n}$ is our hypothesis. Let $5|n$. Since $n$ is even, we have $n=10m$ and hence $5p\equiv5\pmod{20m}$. This gives $p\equiv 1\pmod{4m}$. Also, 
    we have $p\equiv 1\pmod{n}$, that is, $p\equiv 1\pmod{10m}$. Thus, we must have $p\equiv 1\pmod{20m}$, that is, $p\equiv 1\pmod{2n}$, which is again a contradiction 
    to the hypothesis that $p\not \equiv 1\pmod{2n}$. Hence, if $(\frac{5}{2n})_k$ contains a multiple of $p$ then $(\frac{2n-3}{2n})_k$ must contain a multiple of $p$. 
    \par Now, from \eqref{rev-new-eq-1}, we have
\begin{align}\label{rev-new-eq-2}
\dfrac{(\frac{1}{n}-\frac{x}{n})_k}{(1+\frac{x}{n})_k}&=\dfrac{(\frac{1}{n})_k}{(1)_k}[1+\sum_{i\geq1}a_{k,i}x^i]
\end{align}
for $0\leq k\leq\frac{p-1}{n}$. Therefore,
\begin{align}\label{eq-700}
&{_5}F_4\left[\begin{array}{ccccc}
                \frac{1}{n}, & \frac{2n-3}{2n}, & \frac{1-x}{n}, & \frac{1-y}{n}, & \frac{1-p}{n} \vspace{.08cm}\\
                ~ & \frac{5}{2n}, & 1+\frac{x}{n}, & 1+\frac{y}{n}, & 1+\frac{p}{n}
              \end{array};1
\right]_{\frac{p-1}{n}}\notag\\
&\equiv {_5}F_4\left[\begin{array}{ccccc}
                \frac{1}{n}, & \frac{2n-3}{2n}, & \frac{1-x}{n}, & \frac{1-y}{n}, & \frac{1}{n}\vspace{.08cm} \\
                ~ & \frac{5}{2n}, & 1+\frac{x}{n}, & 1+\frac{y}{n}, & 1
              \end{array};1
\right]_{\frac{p-1}{n}}\pmod{p}\notag\\
&=\sum_{k=0}^{\frac{p-1}{n}}\dfrac{(\frac{2n-3}{2n})_k(\frac{1}{n})_k^4}{(\frac{5}{2n})_k(1)_k^3k!}
\left[1+\sum_{i\geq1}a_{k,i}x^i\right]\left[1+\sum_{i\geq1}a_{k,i}y^i\right]\in\mathbb{Z}_p[[x,y]].
\end{align}
If we put $a=\frac{1}{n}$, $b=\frac{2n-3}{2n}$,
$c=\frac{1}{2n}$, $d=\frac{1-x}{n}$, $e=\frac{1-y}{n}$ and
$f=\frac{1-p}{n}$ in Theorem \ref{thm-4} we obtain
\begin{align}\label{eq-202}
&{_{5}}F_4\left[\begin{array}{ccccc}
                \frac{1}{n}, & \frac{2n-3}{2n}, & \frac{1-x}{n}, & \frac{1-y}{n}, & \frac{1-p}{n} \vspace{.08cm}\\
                ~ & \frac{5}{2n}, & 1+\frac{x}{n}, & 1+\frac{y}{n}, & 1+\frac{p}{n}
              \end{array};1
\right]_{\frac{p-1}{n}}\notag\\
&=\dfrac{(\frac{n+1}{n})_{\frac{p-1}{n}}(\frac{n-1}{n}+\frac{x+y}{n})_{\frac{p-1}{n}}}
{(1+\frac{x}{n})_{\frac{p-1}{n}}(1+\frac{y}{n})_{\frac{p-1}{n}}}
{_4}F_3\left[\begin{array}{cccc}
               \frac{2}{n}, & \frac{1-x}{n}, & \frac{1-y}{n} & \frac{1-p}{n} \vspace{.08cm}\\
               ~ & \frac{2}{n}-\frac{x+y+p}{n}, & \frac{5}{2n}, & \frac{2n+1}{2n}
             \end{array};1
\right]_{\frac{p-1}{n}}.
\end{align}
Here, we note that the rising factorial $(\frac{n+1}{n})_{\frac{p-1}{n}}$ contains a factor multiple of $p$, namely
$(\frac{n+1}{n}+\frac{p-1}{n}-1)=\frac{p}{n}$; and if $n\neq4,~p\not\equiv1\pmod{2n}$, then the rising factorials
$(\frac{5}{2n})_{\frac{p-1}{n}}$ and $(\frac{2n+1}{2n})_{\frac{p-1}{n}}$ do not contain a multiple of $p$. This implies that
\begin{align}\label{eq-701}
&_{5}F_4\left[\begin{array}{ccccc}
                \frac{1}{n}, & \frac{2n-3}{2n}, & \frac{1-x}{n}, & \frac{1-y}{n}, & \frac{1-p}{n} \vspace{.08cm}\\
                ~ & \frac{5}{2n}, & 1+\frac{x}{n}, & 1+\frac{y}{n}, & 1+\frac{p}{n}
              \end{array};1
\right]_{\frac{p-1}{n}}\in p\mathbb{Z}_p[[x,y]].
\end{align}
Now, let $x=\zeta p$ and $y=\zeta^2p$. Then \eqref{eq-700} and Lemma \ref{lemma-swisher} yield
\begin{align}
&_{5}F_4\left[\begin{array}{ccccc}
                \frac{1}{n}, & \frac{2n-3}{2n}, & \frac{1-\zeta p}{n}, & \frac{1-\zeta^2p}{n}, & \frac{1-p}{n} \\
                ~ & \frac{5}{2n}, & 1+\frac{\zeta p}{n}, & 1+\frac{\zeta^2p}{n}, & 1+\frac{p}{n}
              \end{array};1
\right]_{\frac{p-1}{n}}\notag\\
&=A_0+\sum_{i\geq1}A_{3i}p^{3i},
\end{align}
where $\displaystyle A_0=\sum_{k=0}^{\frac{p-1}{n}}\dfrac{(\frac{2n-3}{2n})_k(\frac{1}{n})_k^4}{(\frac{5}{2n})_k(1)_k^3k!}$ and
$\displaystyle A_{3i}=\sum_{k=0}^{\frac{p-1}{n}}\dfrac{(\frac{2n-3}{2n})_k(\frac{1}{n})_k^4}{(\frac{5}{2n})_k(1)_k^3k!}H(a_{k,i})$ and
$H(a_{k,i})$ is an integral polynomial in the $a_{k,i}$, where the second subscripts in each monomial add to $3$.
Now, \eqref{eq-701} implies that $A_0$ and $A_{3i}\in p\mathbb{Z}_p$.
Thus, we obtain
\begin{align}
&_{5}F_4\left[\begin{array}{ccccc}
                \frac{1}{n}, & \frac{2n-3}{2n}, & \frac{1-\zeta p}{n}, & \frac{1-\zeta^2p}{n}, & \frac{1-p}{n} \\
                ~ & \frac{5}{2n}, & 1+\frac{\zeta p}{n}, & 1+\frac{\zeta^2p}{n}, & 1+\frac{p}{n}
              \end{array};1
\right]_{\frac{p-1}{n}}\notag\\
&\equiv
\sum_{k=0}^{\frac{p-1}{n}}\dfrac{(\frac{2n-3}{2n})_k(\frac{1}{n})_k^4}{(\frac{5}{2n})_k(1)_k^3k!}\notag\\
&\equiv{_5}F_{4}\left[\begin{array}{ccccc}
                 \frac{1}{n}, & \frac{1}{n}, & \frac{1}{n}, & \frac{1}{n}, & \frac{2n-3}{2n}\vspace{.08cm} \\
                 ~ & \frac{5}{2n}, & 1, & 1, & 1
               \end{array};1
\right]_{\frac{p-1}{n}}\pmod{p^4},\notag
\end{align}
which completes the proof of \eqref{eq-400}.
\par Again, from the right hand side of \eqref{eq-200} we have
\begin{align}\label{eq-204}
&{_5}F_{4}\left[\begin{array}{ccccc}
                 \frac{1}{n}, & \frac{1}{n}, & \frac{1}{n}, & \frac{1}{n}, & \frac{2n-3}{2n}\vspace{.08cm} \\
                 ~ & \frac{5}{2n}, & 1, & 1, & 1
               \end{array};1
\right]_{\frac{p-1}{n}}\notag\\
&\equiv\dfrac{(\frac{n+1}{n})_{\frac{p-1}{n}}(\frac{2}{n})_{\frac{p-1}{n}}(\frac{3}{2n}+
\frac{\zeta p}{n})_{\frac{p-1}{n}}(\frac{2n-1}{2n}+\frac{\zeta p}{n})_{\frac{p-1}{n}}}
{(\frac{5}{2n})_{\frac{p-1}{n}}(\frac{2n+1}{2n})_{\frac{p-1}{n}}
(1+\frac{\zeta p}{n})_{\frac{p-1}{n}}(\frac{1}{n}+\frac{\zeta p}{n})_{\frac{p-1}{n}}}\pmod {p^4}.
\end{align}
Here, we observe that the rising factorial $(\frac{n+1}{n})_{\frac{p-1}{n}}$ contains a factor multiple of $p$, namely
$(\frac{n+1}{n}+\frac{p-1}{n}-1)=\frac{p}{n}$; and if $n\neq4,p\not\equiv1\pmod{2n}$, then no other terms on the right
side of \eqref{eq-204} contain a multiple of $p$.
Proceeding similarly as shown in the proof of Theorem \ref{MT-4}, we deduce that the right side of \eqref{eq-204} becomes
\begin{align}\label{eq-953}
&\dfrac{(\frac{n+1}{n})_{\frac{p-1}{n}}(\frac{2}{n})_{\frac{p-1}{n}}(\frac{3}{2n}+
\frac{\zeta p}{n})_{\frac{p-1}{n}}(\frac{2n-1}{2n}+\frac{\zeta p}{n})_{\frac{p-1}{n}}}
{(\frac{5}{2n})_{\frac{p-1}{n}}(\frac{2n+1}{2n})_{\frac{p-1}{n}}
(1+\frac{\zeta p}{n})_{\frac{p-1}{n}}(\frac{1}{n}+\frac{\zeta p}{n})_{\frac{p-1}{n}}}\notag\\
&=\dfrac{p}{2n}\dfrac{\Gamma_p(1+\frac{p}{n})\Gamma_p(\frac{1}{n}+\frac{p}{n})\Gamma_p(\frac{1}{2n}+\frac{(1+\zeta)p}{n})
\Gamma_p(\frac{2n-3}{2n}+\frac{(1+\zeta)p}{n})}{\Gamma_p(\frac{1}{n})\Gamma_p(\frac{2}{n})
\Gamma_p(\frac{3}{2n}+\frac{\zeta p}{n})\Gamma_p(\frac{2n-1}{2n}+\frac{\zeta p}{n})}\notag\\
&\times\dfrac{\Gamma_p(\frac{5}{2n})\Gamma_p(\frac{1}{2n})\Gamma_p(1+\frac{\zeta p}{n})
\Gamma_p(\frac{1}{n}+\frac{\zeta p}{n})}
{\Gamma_p(\frac{3}{2n}+\frac{p}{n})\Gamma_p(\frac{2n-1}{2n}+\frac{p}{n})
\Gamma_p(\frac{n-1}{n}+\frac{(1+\zeta)p}{n})\Gamma_p(\frac{(1+\zeta)p}{n})}.
\end{align}
Now, arranging the terms with respect to the symmetry of cubic roots of unity, and then applying part (2) of
Proposition \ref{prop-2} and Theorem \ref{thm-long-14} we obtain
\begin{align}\label{eq-1000}
\frac{\Gamma_p(1+\frac{p}{n})\Gamma_p(1+\frac{\zeta p}{n})}{\Gamma_p(-\frac{\zeta^2p}{n})}&=
1+O(p^3),
\end{align}
\begin{align}\label{eq-950}
\dfrac{\Gamma_p(\frac{1}{n}+\frac{p}{n})\Gamma_p(\frac{1}{n}+\frac{\zeta p}{n})}
{\Gamma_p(\frac{n-1}{n}-\frac{\zeta^2p}{n})}&=(-1)^{a_0(1-\frac{1}{n})}
\Gamma_p\left(\frac{1}{n}+\frac{p}{n}\right)\Gamma_p\left(\frac{1}{n}+\frac{\zeta p}{n}\right)
\Gamma_p\left(\frac{1}{n}+\frac{\zeta^2 p}{n}\right)\notag\\
&=(-1)^{1+\frac{p-1}{n}}\Gamma_p\left(\frac{1}{n}\right)^3(1+O(p^3)),
\end{align}
\begin{align}\label{eq-951}
\dfrac{\Gamma_p(\frac{1}{2n}-\frac{\zeta^2p}{n})}{\Gamma_p(\frac{2n-1}{2n}+
\frac{\zeta p}{n})\Gamma_p(\frac{2n-1}{2n}+\frac{p}{n})}
=\Gamma_p\left(\frac{1}{2n}\right)^3\left(1+O(p^3)\right),
\end{align}
and
\begin{align}\label{eq-952}
\dfrac{\Gamma_p(\frac{2n-3}{2n}-\frac{\zeta^2p}{n})}
{\Gamma_p(\frac{3}{2n}+\frac{\zeta p}{n})\Gamma_p(\frac{3}{2n}+\frac{p}{n})}=\Gamma_p\left(\frac{2n-3}{2n}\right)^3
\left(1+O(p^3)\right).
\end{align}
Substituting \eqref{eq-1000}, \eqref{eq-950}, \eqref{eq-951} and \eqref{eq-952} into \eqref{eq-953} we conclude that
\begin{align}\label{eq-954}
&\dfrac{(\frac{n+1}{n})_{\frac{p-1}{n}}(\frac{2}{n})_{\frac{p-1}{n}}(\frac{3}{2n}+
\frac{\zeta p}{n})_{\frac{p-1}{n}}(\frac{2n-1}{2n}+\frac{\zeta p}{n})_{\frac{p-1}{n}}}
{(\frac{5}{2n})_{\frac{p-1}{n}}(\frac{2n+1}{2n})_{\frac{p-1}{n}}
(1+\frac{\zeta p}{n})_{\frac{p-1}{n}}(\frac{1}{n}+\frac{\zeta p}{n})_{\frac{p-1}{n}}}\notag\\
&\equiv(-1)^{1+\frac{p-1}{n}}\frac{p}{2n}\frac{\Gamma_p\left(\frac{1}{n}\right)^2\Gamma_p\left(\frac{1}{2n}\right)^4
\Gamma_p\left(\frac{2n-3}{2n}\right)^3
\Gamma_p\left(\frac{5}{2n}\right)}{\Gamma_p\left(\frac{2}{n}\right)}\pmod {p^4}.
\end{align}
Finally, combining \eqref{eq-400}, \eqref{eq-204} and \eqref{eq-954} we complete the proof of the theorem.
\end{proof}
\begin{proof}[Proof of Theorem \ref{MT-2}]
Let $a=\frac{1}{8}$, $b=\frac{1-\zeta p}{4}$, $c=\frac{1-\zeta^2p}{4}$, $d=\frac{1-\zeta^3p}{4}$, $e=\frac{1-\zeta^4p}{4}$,
$f=\frac{1-p}{4}$, where $\zeta$ is a primitive $5$th root of unity. Substituting all these values in Theorem \ref{thm-2} we have
\begin{align}\label{eq-18}
&{_7}F_6\left[\begin{array}{ccccccc}
               \frac{1}{8}, & \frac{17}{16}, & \frac{1-\zeta p}{4}, & \frac{1-\zeta^2p}{4},
                & \frac{1-\zeta^3p}{4}, & \frac{1-\zeta^4p}{4}, & \frac{1-p}{4} \\
               ~ & \frac{1}{16}, & \frac{7}{8}+\frac{\zeta p}{4}, & \frac{7}{8}+\frac{\zeta^2p}{4},
               & \frac{7}{8}+\frac{\zeta^3p}{4}, & \frac{7}{8}+\frac{\zeta^4p}{4}, & \frac{7}{8}+\frac{p}{4}
             \end{array};1
\right]_{\frac{p-1}{4}}\notag\\
&=\dfrac{(\frac{9}{8})_{\frac{p-1}{4}}(\frac{5+2\zeta p+2\zeta^2p}{8})_{\frac{p-1}{4}}
(\frac{5+2\zeta p+2\zeta^3p}{8})_{\frac{p-1}{4}}(\frac{5+2\zeta^2 p+2\zeta^3p}{8})_{\frac{p-1}{4}}}
{(\frac{7}{8}+\frac{\zeta p}{4})_{\frac{p-1}{4}}(\frac{7}{8}+\frac{\zeta^2 p}{4})_{\frac{p-1}{4}}
(\frac{7}{8}+\frac{\zeta^3 p}{4})_{\frac{p-1}{4}}(\frac{3+2\zeta p+2\zeta^2p+2\zeta^3p}{8})_{\frac{p-1}{4}}}.
\end{align}
Our claim is that
\begin{align}\label{eq-300}
&{_7}F_6\left[\begin{array}{ccccccc}
               \frac{1}{8}, & \frac{17}{16}, & \frac{1-\zeta p}{4}, & \frac{1-\zeta^2p}{4},
                & \frac{1-\zeta^3p}{4}, & \frac{1-\zeta^4p}{4}, & \frac{1-p}{4} \vspace{.08cm}\\
               ~ & \frac{1}{16}, & \frac{7}{8}+\frac{\zeta p}{4}, & \frac{7}{8}+\frac{\zeta^2p}{4},
               & \frac{7}{8}+\frac{\zeta^3p}{4}, & \frac{7}{8}+\frac{\zeta^4p}{4}, & \frac{7}{8}+\frac{p}{4}
             \end{array};1
\right]_{\frac{p-1}{4}}\notag\\
&\equiv{_7}F_6\left[\begin{array}{ccccccc}
                      \frac{1}{8}, & \frac{17}{16}, & \frac{1}{4}, & \frac{1}{4}, & \frac{1}{4}, & \frac{1}{4}, & \frac{1}{4} \vspace{.1cm}\\
                      ~ & \frac{1}{16}, & \frac{7}{8}, & \frac{7}{8}, & \frac{7}{8}, & \frac{7}{8}, & \frac{7}{8}
                    \end{array};1
\right]_{\frac{p-1}{4}}\notag\\&\equiv{_7}F_6\left[\begin{array}{ccccccc}
                      \frac{1}{8}, & \frac{17}{16}, & \frac{1}{4}, & \frac{1}{4}, & \frac{1}{4}, & \frac{1}{4}, & \frac{1}{4} \vspace{.1cm}\\
                      ~ & \frac{1}{16}, & \frac{7}{8}, & \frac{7}{8}, & \frac{7}{8}, & \frac{7}{8}, & \frac{7}{8}
                    \end{array};1
\right]_{\frac{7(p-1)}{8}}\pmod{p^6}.
\end{align}
The last congruence follows easily from the fact that the rising factorials $(\frac{1}{8})_k$ and $(\frac{1}{4})_k$ contain a multiple of $p$
for each $k$ satisfying $\frac{p-1}{4}<k\leq \frac{7(p-1)}{8}$.
Again,
\begin{align}\label{eq-301}
\frac{(\frac{1}{4}-x)_k}{(\frac{7}{8}+x)_k}&=\prod_{j=0}^{k-1}\frac{(\frac{1}{4}+j-x)}{(\frac{7}{8}+j+x)}
=\prod_{j=0}^{k-1}\frac{(\frac{1}{4}+j)(1-\frac{4x}{1+4j})}{(\frac{7}{8}+j)(1+\frac{8x}{7+8j})}\notag\\&
=\frac{(\frac{1}{4})_k}{(\frac{7}{8})_k}[1+a_{k,1}x+a_{k,2}x^2+\cdots]
\end{align}
for some constants $a_{k,1}, a_{k,2}, \ldots$. We observe that the terms $\frac{1}{4}+j$ and $\frac{7}{8}+j$ do not
contain a multiple of $p$ for each $k$
in the range $0\leq k\leq\frac{p-1}{4}$, and hence $a_{k,1}, a_{k,2}, \ldots \in\mathbb{Z}_p$.
Now, \eqref{eq-301} yields
\begin{align}\label{eq-306}
&{_7}F_6\left[\begin{array}{ccccccc}
                      \frac{1}{8}, & \frac{17}{16}, & \frac{1}{4}-x, & \frac{1}{4}-y, & \frac{1}{4}-z, & \frac{1}{4}-w, & \frac{1-p}{4}\vspace{.08cm} \\
                     ~ & \frac{1}{16}, & \frac{7}{8}+x, & \frac{7}{8}+y, & \frac{7}{8}+z, & \frac{7}{8}+w, & \frac{7}{8}+\frac{p}{4}
                    \end{array};1\right]_{\frac{p-1}{4}}\notag\\
                    &\equiv{_7}F_6\left[\begin{array}{ccccccc}
                      \frac{1}{8}, & \frac{17}{16}, & \frac{1}{4}-x, & \frac{1}{4}-y, & \frac{1}{4}-z, & \frac{1}{4}-w, & \frac{1}{4}\vspace{.08cm} \\
                     ~ & \frac{1}{16}, & \frac{7}{8}+x, & \frac{7}{8}+y, & \frac{7}{8}+z, & \frac{7}{8}+w, & \frac{7}{8}
                    \end{array};1\right]_{\frac{p-1}{4}}\pmod{p}\notag\\
                    &=\sum_{k=0}^{\frac{p-1}{4}}\frac{(16k+1)(\frac{1}{8})_k(\frac{1}{4})_k^5}{(\frac{7}{8})_k^5k!}\left[1+\sum_{i\geq1}
                    a_{k,i}x^i\right]\left[1+\sum_{i\geq1}a_{k,i}y^i\right]\left[1+\sum_{i\geq1}a_{k,i}z^i\right]\notag\\
                   &\hspace{.5cm} \times\left[1+\sum_{i\geq1}a_{k,i}w^i\right]
                    \in\mathbb{Z}_p[[x,y,z,w]].
\end{align}
Using Theorem \ref{thm-4} we obtain
\begin{align}
&{_7}F_6\left[\begin{array}{ccccccc}
                      \frac{1}{8}, & \frac{17}{16}, & \frac{1}{4}-x, & \frac{1}{4}-y, & \frac{1}{4}-z, & \frac{1}{4}-w, & \frac{1-p}{4} \vspace{.08cm}\\
                     ~ & \frac{1}{16}, & \frac{7}{8}+x, & \frac{7}{8}+y, & \frac{7}{8}+z, & \frac{7}{8}+w, & \frac{7}{8}+\frac{p}{4}
                    \end{array};1\right]_{\frac{p-1}{4}}\notag\\
&=\frac{(\frac{9}{8})_{\frac{p-1}{4}}(\frac{5}{8}+z+w)_{\frac{p-1}{4}}}
{(\frac{7}{8}+z)_{\frac{p-1}{4}}(\frac{7}{8}+w)_{\frac{p-1}{4}}}\notag\\
&\times {_4}F_3\left[\begin{array}{cccc}
               \frac{5}{8}+x+y, & \frac{1}{4}-z, & \frac{1}{4}-w, & \frac{1-p}{4} \vspace{.08cm}\\
               ~ & \frac{7}{8}+x, & \frac{7}{8}+y, & \frac{5}{8}-z-w-\frac{p}{4}
             \end{array};1
\right]_{\frac{p-1}{4}}.\notag
\end{align}
The rising factorial $(\frac{9}{8})_{\frac{p-1}{4}}$ contains a multiple of $p$, namely $(\frac{9}{8}+\frac{p-1}{8}-1)=\frac{p}{8}$;
and the terms in the denominator are units. Therefore, we have
\begin{align}\label{eq-302}
&{_7}F_6\left[\begin{array}{ccccccc}
                      \frac{1}{8}, & \frac{17}{16}, & \frac{1}{4}-x, & \frac{1}{4}-y, & \frac{1}{4}-z, & \frac{1}{4}-w, & \frac{1-p}{4} \vspace{.1cm}\\
                     ~ & \frac{1}{16}, & \frac{7}{8}+x, & \frac{7}{8}+y, & \frac{7}{8}+z, & \frac{7}{8}+w, & \frac{7}{8}+\frac{p}{4}
                    \end{array};1\right]_{\frac{p-1}{4}}\in p\mathbb{Z}_p[[x,y,z,w]].
\end{align}
Let $x=\frac{\zeta p}{4}, y=\frac{\zeta^2p}{4}, z=\frac{\zeta^3p}{4}$ and $w=\frac{\zeta^4p}{4}$.
Then \eqref{eq-306} and Lemma \ref{lemma-swisher}  yield
\begin{align}
&{_7}F_6\left[\begin{array}{ccccccc}
               \frac{1}{8}, & \frac{17}{16}, & \frac{1}{4}-\frac{\zeta p}{4}, & \frac{1}{4}-\frac{\zeta^2p}{4}, & \frac{1}{4}-\frac{\zeta^3p}{4}, &
               \frac{1}{4}-\frac{\zeta^4p}{4}, & \frac{1}{4}-\frac{p}{4} \vspace{.1cm}\\
               ~ & \frac{1}{16}, & \frac{7}{8}+\frac{\zeta p}{4}, & \frac{7}{8}+\frac{\zeta^2p}{4}, & \frac{7}{8}+\frac{\zeta^3p}{4}, & \frac{7}{8}+\frac{\zeta^4p}{4},
               & \frac{7}{8}+\frac{p}{4}
             \end{array};1
\right]_{\frac{p-1}{4}}\notag\\&=A_0+\sum_{i\geq1}A_{5i}p^{5i},\notag
\end{align}
where $\displaystyle A_0=\sum_{k=0}^{\frac{p-1}{4}}\frac{(16k+1)(\frac{1}{8})_k(\frac{1}{4})_k^5}{(\frac{7}{8})_k^5k!}$ and
$\displaystyle A_{5i}=\sum_{k=0}^{\frac{p-1}{4}}\frac{(16k+1)(\frac{1}{8})_k(\frac{1}{4})_k^5}{4^5(\frac{7}{8})_k^5k!}H(a_{k,i})$;
and $H(a_{k,i})$ is an integral polynomial in the $a_{k,i}$ where the second subscripts in each monomial add to $5$.
Now, \eqref{eq-302} implies that $A_0$ and $A_{5i}$ are in $p\mathbb{Z}_p$.
Therefore, we conclude that
\begin{align}\label{eq-305}
&{_7}F_6\left[\begin{array}{ccccccc}
               \frac{1}{8}, & \frac{17}{16}, & \frac{1-\zeta p}{4}, & \frac{1-\zeta^2p}{4},
                & \frac{1-\zeta^3p}{4}, & \frac{1-\zeta^4p}{4}, & \frac{1-p}{4} \vspace{.1cm}\\
               ~ & \frac{1}{16}, & \frac{7}{8}+\frac{\zeta p}{4}, & \frac{7}{8}+\frac{\zeta^2p}{4},
               & \frac{7}{8}+\frac{\zeta^3p}{4}, & \frac{7}{8}+\frac{\zeta^4p}{4}, & \frac{7}{8}+\frac{p}{4}
             \end{array};1
\right]_{\frac{p-1}{4}}\notag\\
&\equiv\sum_{k=0}^{\frac{p-1}{4}}\frac{(16k+1)(\frac{1}{8})_k(\frac{1}{4})_k^5}{(\frac{7}{8})_k^5k!}\notag\\
&={_7}F_6\left[\begin{array}{ccccccc}
                      \frac{1}{8}, & \frac{17}{16}, & \frac{1}{4}, & \frac{1}{4}, & \frac{1}{4}, & \frac{1}{4}, & \frac{1}{4} \vspace{.08cm}\\
                      ~ & \frac{1}{16}, & \frac{7}{8}, & \frac{7}{8}, & \frac{7}{8}, & \frac{7}{8}, & \frac{7}{8}
                    \end{array};1
\right]_{\frac{p-1}{4}}\pmod{p^6}.
\end{align}
This proves our claim \eqref{eq-300}.
\par If we apply \eqref{new-eq-2} and \eqref{new-eq-3}, and then part (2) of Proposition \ref{prop-1} on the right hand side of
\eqref{eq-18} as in Theorem \ref{MT-4} and Theorem \ref{MT-6}, we deduce that
\begin{align}\label{new-eq-5}
&\dfrac{(\frac{9}{8})_{\frac{p-1}{4}}(\frac{5+2\zeta p+2\zeta^2p}{8})_{\frac{p-1}{4}}
(\frac{5+2\zeta p+2\zeta^3p}{8})_{\frac{p-1}{4}}(\frac{5+2\zeta^2 p+2\zeta^3p}{8})_{\frac{p-1}{4}}}
{(\frac{7}{8}+\frac{\zeta p}{4})_{\frac{p-1}{4}}(\frac{7}{8}+\frac{\zeta^2 p}{4})_{\frac{p-1}{4}}
(\frac{7}{8}+\frac{\zeta^3 p}{4})_{\frac{p-1}{4}}(\frac{3+2\zeta p+2\zeta^2p+2\zeta^3p}{8})_{\frac{p-1}{4}}}\notag\\
&=\dfrac{-p\Gamma_p(\frac{7+2p}{8})\Gamma_p(\frac{3+2p+2\zeta p+2\zeta^2p}{8})
\Gamma_p(\frac{3+2p+2\zeta p+2\zeta^3p}{8})\Gamma_p(\frac{3+2p+2\zeta^2 p+2\zeta^3p}{8})}
{\Gamma_p(\frac{1}{8})\Gamma_p(\frac{5+2\zeta p+2\zeta^2p}{8})\Gamma_p(\frac{5+2\zeta p+2\zeta^3p}{8})
\Gamma_p(\frac{5+2\zeta^2 p+2\zeta^3p}{8})}\notag\\
&\times\dfrac{\Gamma_p(\frac{7}{8}+\frac{\zeta p}{4})\Gamma_p(\frac{7}{8}+\frac{\zeta^2 p}{4})
\Gamma_p(\frac{7}{8}+\frac{\zeta^3 p}{4})\Gamma_p(\frac{3+2\zeta p+2\zeta^2p+2\zeta^3p}{8})}{\Gamma_p(\frac{5}{8}+\frac{(1+\zeta) p}{4})
\Gamma_p(\frac{5}{8}+\frac{(1+\zeta^2) p}{4})\Gamma_p(\frac{5}{8}+\frac{(1+\zeta^3)p}{4})
\Gamma_p(\frac{1+2p+2\zeta p+2\zeta^2p+2\zeta^3p}{8})}.
\end{align}
Rearranging the terms on the right side of \eqref{new-eq-5} with respect to the symmetry of the 5th root of unity, we obtain
\begin{align}\label{eq-900}
&\dfrac{(\frac{9}{8})_{\frac{p-1}{4}}(\frac{5+2\zeta p+2\zeta^2p}{8})_{\frac{p-1}{4}}
(\frac{5+2\zeta p+2\zeta^3p}{8})_{\frac{p-1}{4}}(\frac{5+2\zeta^2 p+2\zeta^3p}{8})_{\frac{p-1}{4}}}
{(\frac{7}{8}+\frac{\zeta p}{4})_{\frac{p-1}{4}}(\frac{7}{8}+\frac{\zeta^2 p}{4})_{\frac{p-1}{4}}
(\frac{7}{8}+\frac{\zeta^3 p}{4})_{\frac{p-1}{4}}(\frac{3+2\zeta p+2\zeta^2p+2\zeta^3p}{8})_{\frac{p-1}{4}}}\notag\\
&=\dfrac{-p\Gamma_p(\frac{7+2p}{8})\Gamma_p(\frac{7}{8}+\frac{\zeta p}{4})\Gamma_p(\frac{7}{8}+\frac{\zeta^2 p}{4})
\Gamma_p(\frac{7}{8}+\frac{\zeta^3 p}{4})}{\Gamma_p(\frac{1}{8})\Gamma_p(\frac{1+2p+2\zeta p+2\zeta^2p+2\zeta^3p}{8})}\notag\\
&\times\dfrac{\Gamma_p(\frac{3+2p+2\zeta p+2\zeta^2p}{8})
\Gamma_p(\frac{3+2p+2\zeta p+2\zeta^3p}{8})\Gamma_p(\frac{3+2p+2\zeta^2 p+2\zeta^3p}{8})}
{\Gamma_p(\frac{5+2\zeta p+2\zeta^2p}{8})\Gamma_p(\frac{5+2\zeta p+2\zeta^3p}{8})
\Gamma_p(\frac{5+2\zeta^2 p+2\zeta^3p}{8})}\notag\\
&\times\dfrac{\Gamma_p(\frac{3+2\zeta p+2\zeta^2p+2\zeta^3p}{8})}{\Gamma_p(\frac{5}{8}+\frac{(1+\zeta) p}{4})
\Gamma_p(\frac{5}{8}+\frac{(1+\zeta^2) p}{4})\Gamma_p(\frac{5}{8}+\frac{(1+\zeta^3)p}{4})}.
\end{align}
Now, using part (2) of Proposition \ref{prop-2} and
Theorem \ref{thm-long-14} on the right side of \eqref{eq-900}, we deduce that
\begin{align}\label{eq-901}
\dfrac{\Gamma_p(\frac{7+2p}{8})\Gamma_p(\frac{7}{8}+\frac{\zeta p}{4})\Gamma_p(\frac{7}{8}+\frac{\zeta^2 p}{4})
\Gamma_p(\frac{7}{8}+\frac{\zeta^3 p}{4})}{\Gamma_p(\frac{1+2p+2\zeta p+2\zeta^2p+2\zeta^3p}{8})}=(-1)^{1+\frac{p-1}{8}}
\Gamma_p\left(\frac{7}{8}\right)^5[1+O(p^5)],
\end{align}
and
\begin{align}\label{eq-902}
&\dfrac{\Gamma_p(\frac{3}{8}+\frac{(1+\zeta+\zeta^2)p}{4})
\Gamma_p(\frac{3}{8}+\frac{(1+\zeta+\zeta^3)p}{4})\Gamma_p(\frac{3}{8}+\frac{(1+\zeta^2+\zeta^3)p}{4})}
{\Gamma_p(\frac{5}{8}+\frac{(\zeta+\zeta^2)p}{4})\Gamma_p(\frac{5}{8}+\frac{(\zeta+\zeta^3)p}{4})
\Gamma_p(\frac{5}{8}+\frac{(\zeta^2+\zeta^3)p}{4})\Gamma_p(\frac{5}{8}+\frac{(1+\zeta) p}{4})}\notag\\
&\times\dfrac{\Gamma_p(\frac{3}{8}+\frac{(\zeta+\zeta^2+\zeta^3)p}{4})}
{\Gamma_p(\frac{5}{8}+\frac{(1+\zeta^2) p}{4})\Gamma_p(\frac{5}{8}+\frac{(1+\zeta^3)p}{4})}\notag\\
&=\Gamma_p\left(\frac{3}{8}\right)^{10}[1+O(p^5)].
\end{align}
Again, using \eqref{eq-901}, \eqref{eq-902} and $\Gamma_p\left(\frac{1}{8}\right)\Gamma_p\left(\frac{7}{8}\right)=(-1)^{1+\frac{p-1}{8}}$,
we obtain
\begin{align}\label{eq-23}
&\dfrac{(\frac{9}{8})_{\frac{p-1}{4}}(\frac{5+2\zeta p+2\zeta^2p}{8})_{\frac{p-1}{4}}
(\frac{5+2\zeta p+2\zeta^3p}{8})_{\frac{p-1}{4}}(\frac{5+2\zeta^2 p+2\zeta^3p}{8})_{\frac{p-1}{4}}}
{(\frac{7}{8}+\frac{\zeta p}{4})_{\frac{p-1}{4}}(\frac{7}{8}+\frac{\zeta^2 p}{4})_{\frac{p-1}{4}}
(\frac{7}{8}+\frac{\zeta^3 p}{4})_{\frac{p-1}{4}}(\frac{3+2\zeta p+2\zeta^2p+2\zeta^3p}{8})_{\frac{p-1}{4}}}\notag\\
&=-p\Gamma_p\left(\frac{7}{8}\right)^6\Gamma_p\left(\frac{3}{8}\right)^{10}\left(1+O(p^5)\right).
\end{align}
Finally, from \eqref{eq-18}, \eqref{eq-305}, \eqref{eq-300} and \eqref{eq-23} we obtain
\begin{align}
{_7}F_6\left[\begin{array}{ccccccc}
               \frac{1}{8}, & \frac{17}{16}, & \frac{1}{4}, & \frac{1}{4},
                & \frac{1}{4}, & \frac{1}{4}, & \frac{1}{4} \vspace{.1cm}\\
               ~ & \frac{1}{16}, & \frac{7}{8}, & \frac{7}{8},
               & \frac{7}{8}, & \frac{7}{8}, & \frac{7}{8}
             \end{array};1
\right]_{\frac{7(p-1)}{8}}\equiv-p\Gamma_p\left(\frac{7}{8}\right)^6\Gamma_p\left(\frac{3}{8}\right)^{10}\pmod{p^6}.\notag
\end{align}
This completes the proof of the theorem.
\end{proof}


\bibliographystyle{amsplain}

\end{document}